\documentclass[review,onefignum,onetabnum]{siamart190516}

\newsiamremark{example}{Example}

\usepackage{lipsum}
\usepackage{amsfonts}
\usepackage{graphicx}
\usepackage{epstopdf}
\usepackage{algorithmic}
\ifpdf
  \DeclareGraphicsExtensions{.eps,.pdf,.png,.jpg}
\else
  \DeclareGraphicsExtensions{.eps}
\fi

\newsiamremark{remark}{Remark}
\newsiamremark{hypothesis}{Hypothesis}
\crefname{hypothesis}{Hypothesis}{Hypotheses}
\newsiamthm{claim}{Claim}

\headers{Inverse Scattering by a Random Periodic Structure}{G. Bao, Y. lin and X. Xu}

\title{Inverse Scattering by a Random Periodic Structure\thanks{Submitted to the editors on Feb 23, 2020.
\funding{The work was supported in part by a NSFC Innovative Group Fund
(No.11621101) and the Fundamental Research Funds for the Central Universities.
}}}

\colorlet{siamexlinkcolor}{black}
\author{Gang Bao\thanks{School of Mathematical Sciences, Zhejiang University, Hangzhou 310027, China
  (\email{baog@zju.edu.cn}, \email{linyiwen@zju.edu.cn}, \email{xxu@zju.edu.cn}).}
\and Yiwen Lin\footnotemark[2]
\and Xiang Xu\footnotemark[2]
}

\usepackage{amsopn}

\ifpdf
\hypersetup{
  pdftitle={Inverse Scattering by a Random Periodic Structure},
  pdfauthor={G. Bao and Y. Lin}
}
\fi

\begin{document}

\maketitle

\begin{abstract}
  This paper develops an efficient numerical method for the inverse scattering problem of a time-harmonic plane wave incident on a perfectly reflecting random periodic structure. The method is based on a novel combination of the Monte Carlo technique for sampling the probability space, a continuation method with respect to the wavenumber, and the Karhunen-Lo$\grave{e}$ve expansion of the random structure, which reconstructs key statistical properties of the profile for the unknown random periodic structure from boundary measurements of the scattered fields away from the structure. Numerical results are presented to demonstrate the reliability and efficiency of the proposed method.
\end{abstract}

\begin{keywords}
  Random periodic structure, inverse scattering, Helmholtz equation,  Karhunen-Lo$\grave{e}$ve expansion, Monte Carlo - continuation - uncertainty quantification method.
\end{keywords}

\begin{AMS}
  78A46, 65C30, 65N21  
\end{AMS}

\section{Introduction}
\label{sec:introduction}
Consider the scattering of a time-harmonic electromagnetic plane wave by a periodic structure, all known as a grating in diffractive optics. The scattering theory in periodic structures has many applications, particularly in the design and fabrication of optical elements such as corrective lenses, antireflective interfaces, beam splitters, and sensors \cite{Dobson1993, Elschner1998}. Throughout, it is assumed that the medium is invariant in the $z$-direction. Thus the model problem of the three-dimensional time harmonic Maxwell equations can be reduced to a simpler model problem of the two-dimensional Helmholtz equation. The more complicated biperiodic problem will be considered in a separate work.

The scattering problems in periodic structures have been studied extensively on both mathematical and numerical aspects. We refer to \cite{BaoDobsonCox} and references therein for the mathematical studies of existence and uniqueness of the diffraction grating problems. Numerical methods can be found in \cite{Bruno93, Chen2003AdaptiveFEMPML, Bao2005Adaptive, Nedelec1991} for either an integral equation approach or a variational approach. A comprehensive review can be found in \cite{Petit1980Theory, Bao2001} on diffractive optics technology and its mathematical modeling as well as computational methods.

This paper is concerned with the numerical solution of the inverse problem, which may be described as follows: given the incident field, determine the periodic structure from a measured reflected field a constant distance away from the structure. The inverse problem arises naturally in the study of optimal design problems in diffractive optics, which is to design a grating structure that gives rise to some specified far-field patterns \cite{Dobson1993, Elschner1998, Bao2001optimal}.  A number of numerical methods have been developed to solve these inverse problems.  Ito and Reitich \cite{Ito1997} proposed a high-order perturbation approach based on the methods of variation of boundaries. In \cite{Arens2003}, Arens and Kirsch applied the factorization method to scattering by a periodic surface. Iterative regularization methods were developed by Hettlich in \cite{Hettlich2002} based on the shape derivatives with respect to the variations of the boundary. Bruckner and Elschner \cite{Bruckner2003} gave a two-step optimization algorithm to reconstruct the grating profile. See also Elschner $et~al.$ \cite{EGR} for a reconstruction algorithm of a two-dimensional periodic
structure based on finite elements and optimization techniques.
More recently, Bao $et~al.$ \cite{bao2012ComputationalInverse} presented an efficient continuation method to capture both the macro and micro structures of the grating profiles with multiple frequency data. The method was further extended to the case of phaseless data in \cite{Bao2013phaseless}.  Related uniqueness results and stability estimates for the inverse diffraction problem were obtained in \cite{Kirsch1994, Bao1994, Bao1995, Bao2011unique}.

Existing studies of the inverse problem mostly assume that the periodic structure is deterministic and only the noise level of the measured data is considered for the inverse problem. In practice, however, there is a level of uncertainty of the scattering surface, e.g. the grating structure may have manufacturing defects or it may suffer other possible damages from regular usage. Therefore, in addition to the noise level of measurements, the random surface itself also influences the measured scattered fields. Considering the inverse diffraction problem by random surfaces is closer to reality. In fact, surface roughness measurements are of great significance for the functional performance evaluation of machined parts and design of micro-optical elements. In this scenario, the analytical methods are insufficient any more to deliver desired error tolerances. Little is known in mathematics or computation about solving inverse problems of determining random surfaces. One challenge lies in the fact that the scattered fields depend nonlinearly on the surface, which makes the random surface reconstruction problem extremely difficult. Another challenge is to understand to what extend the reconstruction could be made. In other words, what statistical quantities of the profile could be recovered from the measured data? Until now, only some initial progress has been made. In \cite{bao2014near}, an asymptotic technique has been applied to obtain the analytical solutions with reasonable accuracy in the case of a small perturbation of a plane surface. B-splines have been employed in \cite{firoozabadi2006subsurface} to recover the unknown rough surface, which is represented by the control points for the spline function. The Monte Carlo method in \cite{johnson1996backscattering} has been adopted to calculate the statistics of the scattered waves. Shi $et~al.$ \cite{SHI2017} have used an ultrasonic methodology to evaluate the correlation length function of the randomly rough surfaces in solids.

In this work, we propose an efficient numerical method to reconstruct the random periodic structure from mutli-frequency scattered fields away from the structure, in the sense that three critical statistical properties, namely the expectation, root mean square, and correlation length of the random structure may be reconstructed. Our method (MCCUQ method) is based on a novel combination of the Monte Carlo technique for sampling the probability space, a continuation method with respect to the wavenumber, and the Karhunen-Lo$\grave{e}$ve expansion of the random structure. Throughout, the periodic structure is assumed to be a stationary Gaussian process, which is reasonable for many practical situations. For the representation of the random structure, we employ the Karhunen-Lo$\grave{e}$ve expansion in this work , due to the fact that the expansion is optimal in the sense that the mean-square error of the truncation of the expansion after finite terms is minimal \cite{bookUQ}. Also, our work may be viewed as an extension of Bao $et~al.$ \cite{bao2012ComputationalInverse}, where an effective continuation method was first introduced to reconstruct the (deterministic) periodic structure from the multiple frequency data.

It is noted that the Bayesian inversion technique has been developed for solving other inverse problems with randomness \cite{Stuart2017, ZhouBayesian}.
Under the Bayesian framework, the sample size for the Markov Chain Monte Carlo (MCMC) method is required to be sufficiently high, usually at least $10^4$ or even $10^5$. Additionally, the computational cost of solving the nonlinear inverse problems could be very expensive.
Nevertheless, it is difficult for the Bayesian inversion method to derive an explicit expression for the estimation of the critical parameters in the correlation. In this work, we recover the critical parameters by using the properties of the Karhunen-Lo$\grave{e}$ve expansion and its covariance function, and derive an explicit recovery formula to evaluate the two critical parameters, i.e., the root mean square and the correlation length of the random surface, with the Gussian process $f$ and the uniformly distributed noise.

The rest of the paper is organized as follows. The model problem is introduced in \cref{sec:problem}. The Monte Carlo - Continuation - Uncertainty Quantification reconstruction method is presented and discussed in \cref{sec:method}. Several numerical results are presented in \cref{sec:experiments} to demonstrate the accuracy and reliability of the method.  The paper is concluded with some general remarks.

\section{Model problem}
\label{sec:problem}

\begin{figure}[htbp]
	\centering
	\includegraphics[width=13cm]{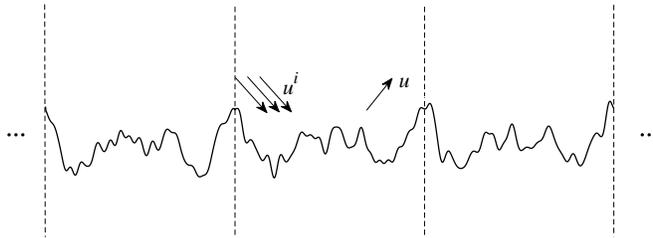}\\
	\caption{Problem geometry}\label{fig:geometry}
\end{figure}

Consider the scattering of an incident field by a periodic random surface
$$\Gamma_f := \{(x,y)\in \mathbb{R}^2|y=f(\omega;x),\omega \in \Omega\},$$
where the rough surface is a stationary Gaussian process as shown in \Cref{fig:geometry}. Here, for a probability space $(\Omega, \mathcal{F}, \mu)$,
$\omega \in \Omega$ denotes the random sample, $(x,y)\in \mathbb{R}^2$ are the spatial variables, and the random surface $f: \Omega \times \mathcal{X}   \rightarrow \mathbb{R}$ is the sum of a deterministic funciton $\tilde{f}(x)$ and a stationary Gaussian process with a continuous and bounded covariance function $c(x-y)$.  The deterministic function is further assumed to be $\Lambda$-periodic, that is,  $\tilde{f}(x+\Lambda)=\tilde{f}(x)$, and the covariance operator $C_f: L^{2}(\mathcal{X}, \mathrm{d} x ; \mathbb{R}) \rightarrow L^{2}(\mathcal{X}, \mathrm{d} x ; \mathbb{R})$ of the random surface $f$ is defined by
$$\left(C_f\varphi\right)(x):= \int_\mathcal{X} c(x-y) \varphi(y) \mbox{d}y,$$
with $\mathcal{X}=[0,\Lambda]. $

For the representation of the random structure, we choose the
Karhunen-Lo$\grave{e}$ve (KL) expansion\cite{feng2018efficientMCTFE,bookUQ}. The KL expansion of a square-integrable stochastic process $f: \Omega \times \mathcal{X}   \rightarrow \mathbb{R}$ is a particularly suitable spectral decomposition for uncertainty quantifications since the mean square error of the truncation of the expansion after finite terms is minimal. In addition, the expansion with respect to the orthonormal basis enjoys many useful properties. Specifically, we have
\begin{equation}
\begin{aligned}
f(\omega;x)&=\tilde{f}(x) + \displaystyle \sum_{j=0}^{\infty}\sqrt{\lambda_j} \xi_{j}(\omega)\varphi_{j}(x),
\end{aligned}
\end{equation}
where $\tilde{f}(x)$ is a $\Lambda$-periodic deterministic function. The orthonormal eigenfunctions of the covariance operator $C_f$ are defined by
\begin{eqnarray}
\varphi_{j}(x)=\left\{\begin{array}{ll}
{\sqrt{\displaystyle\frac{1}{\Lambda}},} & {j=0,} \\ {\sqrt{\displaystyle\frac{2}{\Lambda}} \cos \left(\displaystyle\frac{2 j \pi x}{\Lambda}\right),} & {j>1, \text { even, }} \\ {\sqrt{\displaystyle\frac{2}{\Lambda}} \sin \left(\displaystyle\frac{2 j \pi x}{\Lambda}\right),} & {j>1, \text { odd, }}\end{array}\right.
\end{eqnarray}
where the corresponding eigenvalues $\{ \lambda_j \}_{j=0}^{\infty}$ are arranged in a descending order, and $\{\xi_j \}_{j=0}^{\infty}$ is a random variable with zero mean and unit covariance.
It follows that the KL expansion of the random process $f(\omega;x)$ may be written as
\begin{equation}\label{eq:fKL}
\begin{aligned}
f(\omega;x)&=\tilde{f}(x)+\sqrt{\lambda_0}\xi_0(\omega) \sqrt{\dfrac{1}{\Lambda}} \\
& + \displaystyle \sum_{j=1}^{\infty}\sqrt{\lambda_j}\left( \xi_{j,s}(\omega)\sqrt{\dfrac{2}{\Lambda}} \sin\left(\dfrac{2j\pi x}{\Lambda}\right)
+ \xi_{j,c}(\omega)\sqrt{\dfrac{2}{\Lambda}} \cos\left(\dfrac{2j\pi x}{\Lambda}\right) \right),
\end{aligned}
\end{equation}
where $\xi_0$, $\xi_{j,s}$ and $\xi_{j,c}$ are independent and identically distributed random variables with zero mean and unit covariance.

Since the structure is a Gaussian random surface,  the covariance function \cite{bookRandomSurface, bookcovariance} takes the following form
\begin{equation}
c(x-y)=\sigma^2 exp(-\dfrac{|x-y|^2}{l^2}), 0<l\ll \Lambda,
\end{equation}
where $\sigma$ is the root mean square of the surface, and $l$ is the correlation length.

We next state the model problem. Consider a plane wave incident on a periodic structure $\Gamma_f$ which is characterized by the wavenumbers $\kappa$ ruled on a perfect conductor. It is assumed that the medium and material are invariant in the $z$ direction, and periodic in the $x$ direction with a period $\Lambda$. Above the structure, the medium is assumed to be homogeneous. In this two-dimensional setting, for simplicity, we consider the transverse electric (TE) polarization case.

More precisely, let
\begin{equation}
u^i=e^{i\alpha x-i\beta y}
\end{equation}
be an incident plane wave , where
$ \alpha = \kappa \sin\theta,~\beta=\kappa \cos\theta$,
$\theta\in(-\pi/2,\pi/2)$ is the angle of incidence with respect to the positive $y$-axis. Let the total field
\begin{equation*}
u^{total}(\omega;\cdot)=                        u^i+u(\omega;\cdot), \ \mbox{in}~D_f,
\end{equation*}
where $u$ is the scattered field and $D_f= \{y>f(\omega; x)\}$.
Then, the time-harmonic Maxwell equations can be reduced to the two dimensional Helmholtz equation \cite{BaoDobsonCox}
\begin{equation}\label{eq:vDf}
\Delta u(\omega;\cdot)+\kappa^2 u(\omega;\cdot) = 0, \  \mbox{in}~D_f.
\end{equation}
In addition, as usual, due to the uniqueness consideration, we seek for the quasi-periodic solutions which satisfy \begin{equation}\label{eq:qp}
u(\omega;x+\Lambda,y)= e^{i\alpha \Lambda} u(\omega;x,y),\  \mbox{in}~D_f.
\end{equation}
Moreover, since the medium below $\Gamma_f$ is perfectly electric conducting, from physics, the following surface condition holds
\begin{equation}\label{eq:ic}
u(\omega;\cdot)+u^i=0,\ \mbox{on}\ \Gamma_f.
\end{equation}

For each given sample $\omega$, it follows from the Rayleigh expansion that $u(\omega;\cdot)$ can be written as
\begin{equation}
u(\omega;\cdot)=\sum_{n\in\mathbb{Z}}A_n e^{i\alpha_n x + i \beta_n y},\quad y>\max_{x\in(0,\Lambda)}f(\omega;x),
\end{equation}
where
$$\alpha_n:=\alpha+\dfrac{2\pi n}{\Lambda},~ \beta_n:=\left\{\begin{array}{cc}
\sqrt{\kappa^2-\alpha_n^2}, & |\kappa|>|\alpha_n| \\
i\sqrt{\alpha_n^2-\kappa^2}, & |\kappa|<|\alpha_n|
\end{array}\right..$$
Denote
\begin{equation}
D=\{ (x,y)\in\mathbb{R}^2~|~f(\omega;x)<y<y_0,\forall\omega\in\Omega \},
\end{equation}
with the boundary given by
\begin{equation}
\Gamma_0=\{(x,y)\in\mathbb{R}^2~|~y=y_0\}
\end{equation}
and
$$y_0>\max_{\omega\in\Omega,~x\in(0,\Lambda)}f(\omega;x).$$

Therefore, it follows from Equation \cref{eq:vDf}-\cref{eq:ic} that the random surface scattering model may be stated as the following boundary value problem for the scattered field:
\begin{eqnarray}\label{eq:u}
\left\{\begin{array}{ll}
\Delta u(\omega;\cdot)+\kappa^2 u(\omega;\cdot) = 0, & \mbox{in}~ \Omega \times D,\\
u(\omega;\cdot) + u^i = 0, & \mbox{on}~ \Omega \times \Gamma_f,\\
u(\omega;x+\Lambda,y)=e^{i\alpha \Lambda} u(\omega;x,y), & \mbox{in}~ \Omega \times D.
\end{array}
\right.
\end{eqnarray}

This work is to study the inverse problem of reconstructing a perfectly reflecting random periodic structure.

{\bf Inverse Problem:} For each sample $\omega$, given the incident wave $u^{i}$, determine the KL expansion random structure $y=f(\omega;x)$  as well as the two important statistical properties of the structure,  the root mean square $\sigma$ and the correlation length $l$, from the multiple frequency measurements of the scattered fields $u(\omega;x,y_0)$.

\section{Reconstruction method}
\label{sec:method}
For each given sample $\omega$, we begin with the single-layer potential representation for the scattered field
\begin{equation}
u(\omega;x,y)=\int_{0}^{\Lambda} \phi(\omega;s) G(x,y;s,0) \mbox{d}s,
\end{equation}
where  $\phi\in \mathrm{L}^2(\Omega; \mathrm{L}^2(0,\Lambda))$ is an unknown periodic density function and  $G$ is a quasi-periodic Green function given explicitly as
\begin{equation}\label{Green}
G(x,y;s,t)=\dfrac{i}{2\pi} \sum_{n\in \mathbb{Z}} \dfrac{1}{\beta_n} e^{i\alpha_n (x-s) + i\beta_n |y-t|},\quad (x,y)\neq (s,t).
\end{equation}
Then
\begin{equation}
u(\omega;x,y_0)=\int_{0}^{\Lambda} \phi(\omega;s) G(x,y_0;s,0) \mbox{d}s,
\end{equation}

Since $u$ is quasi-periodic, we can expand
\begin{equation}
u(\omega;x,y_0)=\sum_{n\in \mathbb{Z}} u_n(\omega) e^{i\alpha_n x},
\end{equation}
where
\begin{equation}
u_n(\omega)=\dfrac{1}{\Lambda}\int_0^{\Lambda} u(\omega;x,y_0)e^{-i\alpha_n x}\mbox{d}x.
\end{equation}
On the other hand, the periodic function $\phi$ with periodicity $\Lambda$ can be expanded as
\begin{equation}\label{varphi}
\phi(\omega;s)=\sum_{n\in \mathbb{Z}} \phi_n(\omega) e^{i\alpha_n s}.
\end{equation}
Therefore, we have
\begin{equation}\label{varphin}
\phi_n(\omega)=-i\dfrac{2\pi}{\Lambda}\beta_n u_n(\omega) e^{-i\beta_n y_0}.
\end{equation}

Define the operator $T_f:\mathrm{L}^2(\Omega; \mathrm{L}^2(0,\Lambda)) \rightarrow \mathrm{L}^2(\Omega; \mathrm{L}^2(0,\Lambda))$:
\begin{equation}\label{Tfvarphi}
(T_f\phi)(\omega;x)=\int_0^{\Lambda} \phi(\omega;s) G(x,f(\omega;x);s,0)\mbox{d}s.
\end{equation}
Substituting Equations \cref{varphi} and \cref{varphin} and the quasi-periodic Green function \cref{Green} into Equation \cref{Tfvarphi} gives
\begin{equation}
(T_f\phi)(\omega;x)=\sum_{n\in\mathbb{Z}}\psi_n(\omega) e^{i\alpha_n x + i\beta_n f(\omega;x)}
\end{equation}
where
$$\psi_n(\omega)=u_n(\omega)e^{-i\beta_n y_0}.$$
In practice, in order to restrict the exponential growth of noise, we take the regularization \cite{Heinz1996Regularizaion}  here, that is,
\begin{equation}
\psi_n(\omega)=\left\{\begin{array}{cc}
u_n(\omega)e^{-i\beta_n y_0}, & \mbox{for}\  \kappa>|\alpha_n|, \\
u_n(\omega)\dfrac{e^{i\beta_n y_0}}{e^{2i\beta_n y_0} + \gamma} & \mbox{for}\  \kappa<|\alpha_n|,
\end{array}
\right.
\end{equation}
where $\gamma$ is some positive regularization parameter.

\subsection{The nonlinear problem}
Recalling that $u(\omega;\cdot) + u^i = 0$ in \cref{eq:u}, we may study the nonlinear problem
\begin{equation*}
\left\| (T_f\phi)(\omega;x) + u^i(x,f(\omega,x)) \right\|_{\mathrm{L}^2(\Omega; \mathrm{L}^2(0,\Lambda))}^2=0.
\end{equation*}
Substituting the expansion for the operator $T_f$, we have
\begin{equation*}
\left\|\sum_{n\in\mathbb{Z}}\psi_n(\omega) e^{i\alpha_n x + i\beta_n f(\omega;x)} +e^{i\alpha x-i\beta f(\omega,x)} \right\|_{\mathrm{L}^2(\Omega; \mathrm{L}^2(0,\Lambda))}^2=0.
\end{equation*}
Since the summation decreases exponentially with respect to $|n|$, we choose a sufficiently large $N$, truncate it into finite summation, and consider the numerical solution of the following nonlinear equation
\begin{equation}\label{NonlinearProblem}
\left\|\sum_{n=-N}^{N} \psi_n(\omega) e^{i\alpha_n x + i\beta_n f(\omega;x)} +e^{i\alpha x-i\beta f(\omega,x)} \right\|_{\mathrm{L}^2(\Omega; \mathrm{L}^2(0,\Lambda))}^2=0.
\end{equation}

\subsection{The Monte Carlo technique}
\label{subsec:MC}
In order to approximate the solution of the nonlinear equation  \cref{NonlinearProblem}, the Monte Carlo technique is employed to sample the probability space, which allows us to derive useful statistical properties for the reconstruction.

Let $M$ be a positive integer which denotes the number of realizations. For each sample $\omega_m,~m=1,2,...,M$, we denote the reconstructed random surface by $f_m(x), m=1,2,...,M$. Thus the algorithm yields an approximation of $\mathbb{E}(f)$ given by
\begin{equation}
\bar{f}(x)=\dfrac{1}{M}\sum_{m=1}^M f_m(x).
\end{equation}
Similarly, the standard deviation of the sampled solution can be computed by using the formula
\begin{equation}
s_f(x)=\displaystyle\sqrt{\dfrac{1}{M}\sum_{m=1}^M (f_m(x)-\bar{f})^2},
\end{equation}
the expectation of which is an approximation of the root mean square of the random structure.

Recalling that
\begin{equation*}
\begin{aligned}
f(\omega_m;x)&=\tilde{f}(x)+\sqrt{\lambda_0}\xi_0(\omega_m) \sqrt{\dfrac{1}{\Lambda}} \\
& + \sum_{j=1}^{\infty}\sqrt{\lambda_j}\left( \xi_{j,s}(\omega_m)\sqrt{\dfrac{2}{\Lambda}} \sin\left(\dfrac{2j\pi x}{\Lambda}\right)
+ \xi_{j,c}(\omega_m)\sqrt{\dfrac{2}{\Lambda}} \cos\left(\dfrac{2j\pi x}{\Lambda}\right) \right),
\end{aligned}
\end{equation*}
$\mathbb{E}(\xi_0) = \mathbb{E} (\xi_{j,s}) = \mathbb{E} (\xi_{j,c}) = 0$ when the sample size $M$ is large enough, and the profile $\tilde{f}(x)$ is a real $\Lambda$-period function. Thus the mean value $\bar{f}$ admits a Fourier series expansion. Without loss of generality, we take the period $\Lambda$ to be $2\pi$ from now on. Hence $\bar{f}$ has the following expansion
\begin{equation}\label{eq:fbar}
\bar{f}(x)=\bar{c}_0 + \sum_{p=1}^{\infty} [\bar{c}_{2p-1}\cos(px)+\bar{c}_{2p}\sin(px)],
\end{equation}
where
\begin{equation*}
\bar{c}_p=\dfrac{1}{M}\sum_{m=1}^M c_{m,p},~~ p=0,1,2,...
\end{equation*}
Therefore, the first step of our method is to determine each coefficient $\bar{c}_p~(p=0,1,...)$ to reconstruct the mean structure profile of the random surface.

\subsection{The Monte Carlo - continuation (MCC) method} \label{subsec:MCC}
Because of the nonlinearity of Equation \cref{NonlinearProblem}, we propose a Monte Carlo - continuation (MCC) method to recursively reconstruct these Fourier coefficients for all samples.

It is evident that the finite expansion of Equation \cref{eq:fbar} can reasonably approximate the random periodic grating profile. Here, we choose a prescribed wavenumber not exceeding the wavenumber $\kappa$, and the number of Fourier expansion modes is taken not larger than the prescribed wavenumber.

For each sample $\omega_m,~m=1,2,...,M$, we set the initial approximation $c_{m,0}=y_0$ and $c_{m,p}=0, p=1,2,...,2k_{max}$, where $k_{max}$ is taken to be the largest integer that is smaller than or equal to the prescribed wavenumber.
Denote the vector $$\mathbf{c}_{m,k_{max}}=[c_{m,0},c_{m,1}...,c_{m,2k_{max}}]^{\mathrm{T}}.$$

We first choose an initial value for the wavenumber $\kappa_1\ (\kappa_1<\kappa)$, and seek for an approximation to the profile $f(\omega_m;x)$ by the Fouries series
\begin{equation}
f_{m,k}(x)=c_{m,0}+\sum_{p=1}^{k}\left[c_{m,2 p-1} \cos (p x)+c_{m, 2 p} \sin (p x)\right],
\end{equation}
where $k$ is taken to be the largest integer that is smaller than or equal to the wavenumber $\kappa_1$.
Denote $\mathbf{c}_{m,k}$ the first $2k+1$ terms of $\mathbf{c}_{m,k_{max}}$, i.e.  $\mathbf{c}_{m,k}=[c_{m,0},c_{m,1}...,c_{m,2k}]^{\mathrm{T}}$ with
$c_{m,0}=y_0$ and $c_{m,p}=0, p=1,2,...,2k$.

For the incident angle $\theta_l\in (-\pi/2, \pi/2), l=1,2,...,L$, $\alpha=k sin(\theta_l)$, $\beta=k cos(\theta_l)$, define
$$J_l(\mathbf{c}_{m,k}) =\left\|\displaystyle\sum_{n=-N}^{N} \psi_{n,l}(\omega_m) e^{i\alpha_n x + i\beta_n f(\omega_m;x)} +e^{i\alpha x-i\beta f(\omega_m,x)} \right\|_{\mathrm{L}^2(0,\Lambda)}^2,\   l=1,2,...,L.$$
Denote $\mathbf{J}(\mathbf{c}_{m,k}) = [ J_1(\mathbf{c}_{m,k}), ..., J_L(\mathbf{c}_{m,k}) ]^{\mathrm{T}}$,
then the nonlinear equation \cref{NonlinearProblem} can be reformulated as
\begin{equation}
\mathbf{J}(\mathbf{c}_{m,k})=0,
\end{equation}
where $\mathbf{J}: \mathbb{R}^{2k+1} \rightarrow \mathbb{R}^{L}$.
In order to reduce the computational cost and instability, for each sample, we consider the nonlinear Landweber iteration \cite{Heinz1996Regularizaion}
\begin{equation}
\mathbf{c}_{m,k}^{(t+1)}=\mathbf{c}_{m,k}^{(t)}- \eta_k \mathbf{DJ}^{\mathrm{T}}(\mathbf{c}_{m,k}^{(t)}) \mathbf{J}(\mathbf{c}_{m,k}^{(t)}), \ t=0,1,2,...
\end{equation}
where
$\eta_k$ is a relaxation parameter dependent on the wavenumber, the initial vector $\mathbf{c}_{m,k}^{(0)}=\mathbf{c}_{m,k}$, and
the Jacobi matrix
$$\mathbf{DJ}=\left( \dfrac{\partial J_l}{\partial c_{m,p}} \right)_{l=1,2,....,L,p=0,1,...,2k}$$
can be computed explicitly.
In order to update the reconstructed grating profile function by recursive linearization, we set the resulting coefficients as $\mathbf{c}_{m,k}$. Namely, update the first $2k+1$ terms of $\mathbf{c}_{m,k_{max}}$ by using the above nonlinear Landweber method at the wavenumber $\kappa_1$.

Next, by increasing the wavenumber $\kappa_1$ to $\kappa_2$ ($\kappa_2<\kappa$), we seek for a new approximation to the profile $f(\omega_m;x)$ by the Fourier series
\begin{equation}
f_{m,k^{\prime}}(x)=c_{m,0}+\sum_{p=1}^{k^{\prime}}\left[c_{m,2 p-1} \cos (p x)+c_{m, 2 p} \sin (p x)\right],
\end{equation}
and determine the coefficients $\mathbf{c}_{m,k^{\prime}}$,
where $k^{\prime} \ (k^{\prime}>k)$ is again taken to be the largest integer that is smaller than or equal to the wavenumber $\kappa_2$.
Denote $\mathbf{c}_{m,k^{\prime}} = [c_{m,0},c_{m,1}...,c_{m,2k^{\prime}}]^{\mathrm{T}}$ where
\begin{align*}
c_{m,p}:=\left\{\begin{array}{ll}
{c_{m,p},} & {\text { for }\ 0 \leq p \leq 2k,} \\
{0,} & {\text { for }\ 2k<p\leq 2k^{\prime},}
\end{array}\right.
\end{align*}
and we repeat the above Landweber iteration with the approximation to the profile $f(\omega_m;x)$ at the previous smaller wavenumber as our starting point,  i.e.  $$\mathbf{c}_{m,k'}^{(0)}=\mathbf{c}_{m,k'}.$$

We may repeat the above procedure until the prescribed wavenumber (smaller than $\kappa$) is reached. The resulting coefficients $\mathbf{c}=[\bar{c}_0, \bar{c}_1, ..., \bar{c}_{2k_{max}}]^{\mathrm{T}}$ are deduced by taking the expectation of $\mathbf{c}_{m,p}\ (m=1,2,...,M)$, i.e.
$$\mathbf{c} = \dfrac{1}{M}\displaystyle\sum_{m=1}^M \mathbf{c}_{m,p}.$$
Therefore, the mean grating profile of the random periodic structure is reconstructed as the following finite Fourier series expansion
\begin{equation}\label{eq:fbarfinite}
\bar{f}(x)=\bar{c}_0 + \sum_{p=1}^{k_{max}} [\bar{c}_{2p-1}\cos(px)+\bar{c}_{2p}\sin(px)],
\end{equation}
with
\begin{equation*}
\bar{c}_p=\dfrac{1}{M}\sum_{m=1}^M c_{m,p},~~ p=0,1,...,2k_{max}.
\end{equation*}

Next, we discuss some key practical implementation issues of our MCC method.
On the one hand, we use multi-frequency scattered field data for a range of incident angles to realize the reconstruction of random periodic structure.
At each recursion, the initial value of the high frequency reconstruction depends on the previous lower frequency information.
Our numerical experiments have shown that the present MCC method converges for a larger class of rough surfaces than the usual Newton's method using the same initial value.
On the other hand, we take the relaxation parameter dependent on the wavenumber. In fact our algorithm converges for a larger class of surfaces, provided that the relaxation parameter is chosen to be smaller as the frequency increases.

Our algorithm seems to be convergent and stable. As the result, a good approximation to the random periodic structure at a higher frequency can be obtained at each recursion. The precise description of our method to reconstruct the mean profile $\bar{f}$ is given in \cref{alg:MCC}. It should be pointed that since all of the samples in the probability space are independent and identically distributed, the Monte Carlo sampling process can be performed well in parallel.

\begin{algorithm}
	\caption{MCC method}
	\label{alg:MCC}
	\begin{algorithmic}
		\STATE{Define $\mathbf{c}=[\bar{c}_0, \bar{c}_1, \bar{c}_2,..., \bar{c}_{2k_{max}}]^{\mathrm{T}}$}
		\STATE{Set $\mathbf{c}=0$ (initializing)}
		\FOR{$m=1,2,...,M$ (sampling)}
		\STATE{Generate $\omega_m$, determine $y_0$ }
		\STATE{Define $\mathbf{c}_{m,k_{max}}=[c_{m,0},c_{m,1}...,c_{m,2k_{max}}]^{\mathrm{T}}$}
		\STATE{Set $c_{m,0}=y_0, c_{m,p}=0, p=1,2,...,2k_{max}$ (initializing)}
		\FOR{$k=1,2,...,k_{max}$}
		\STATE{Define $\mathbf{c}_{m,k}^{(0)}=[c_{m,0},c_{m,1}...,c_{m,2k}]^{\mathrm{T}}$}
		\FOR{$t=0,1,2,...T$ (iteration)}
		\STATE{$\alpha=k sin(\theta_l)$,~$\beta=k cos(\theta_l)$,~$\theta_l\in (-\pi/2, \pi/2), l=1,2,...,L$}
		\STATE{Define
			$J_l(\mathbf{c}_{m,k}) =\left\|\displaystyle\sum_{n=-N}^{N} \psi_{n,l}(\omega_m) e^{i\alpha_n x + i\beta_n f(\omega_m;x)} +e^{i\alpha x-i\beta f(\omega_m,x)} \right\|_{\mathrm{L}^2(0,\Lambda)}^2$}
		\STATE{Define $\mathbf{J}(\mathbf{c}_{m,k}) = [ J_1(\mathbf{c}_{m,k}), ..., J_L(\mathbf{c}_{m,k}) ]^{\mathrm{T}}$}
		\STATE{Calculate the Jacobi matrix $\mathbf{DJ}=\left( \dfrac{\partial \mathbf{J}_l}{\partial c_{m,p}} \right)_{l=1,2,....,L,p=0,1,...,2k}$}
		\STATE{$\mathbf{c}_{m,k}^{(t+1)}=\mathbf{c}_{m,k}^{(t)}- \eta_k \mathbf{DJ}^{\mathrm{T}}(\mathbf{c}_{m,k}^{(t)}) \mathbf{J}(\mathbf{c}_{m,k}^{(t)})$}
		\ENDFOR
		\STATE{$c_{m,p}:=c_{m,p}^{(T+1)}, p=0,1,...,2k$}
		\ENDFOR
		\ENDFOR
		\STATE{$\mathbf{c} := \dfrac{1}{M}\displaystyle\sum_{m=1}^M \mathbf{c}_{m,k_{max}}$}
		\RETURN $\mathbf{c}$
	\end{algorithmic}
\end{algorithm}

\subsection{Monte Carlo - continuation - uncertainty quantification (MCCUQ) method}
\label{subsec:MCCUQ}
Based on the mean grating profile $\bar{f}(x)$ reconstructed in \cref{subsec:MCC}, we can realize the reconstruction of the random process $f(\omega;x)$ including the two key statistics characterizing the random structure. Recalling that
\begin{equation*}
\begin{aligned}
f(\omega;x)&=\tilde{f}(x)+\sqrt{\lambda_0}\xi_0(\omega) \sqrt{\dfrac{1}{\Lambda}} \\
& + \displaystyle \sum_{j=1}^{\infty}\sqrt{\lambda_j}\left( \xi_{j,s}(\omega)\sqrt{\dfrac{2}{\Lambda}} \sin\left(\dfrac{2j\pi x}{\Lambda}\right)
+ \xi_{j,c}(\omega)\sqrt{\dfrac{2}{\Lambda}} \cos\left(\dfrac{2j\pi x}{\Lambda}\right) \right),
\end{aligned}
\end{equation*}
we have $f(\omega;x) - \tilde{f}(x)$ is a stationary Gaussian process.

\begin{figure}[ht]
	\centering
	\includegraphics[width=10cm]{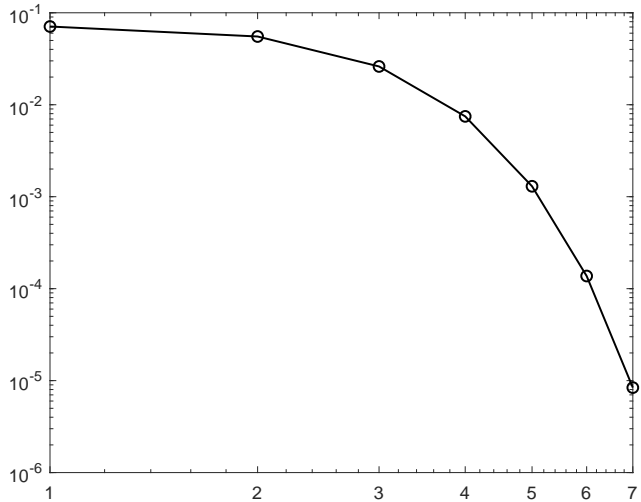}
	\caption{Eigenvalues of the covariance when $\Lambda=2\pi, \sigma=1/5$ and $l=1$.}\label{fig:lambda}
\end{figure}

By simple calculations, we can see that if the eigenvalues $\{\lambda_j\}_{j\in \mathbb{N}}$ of the covariance function are in a descending order, then we only need the first several terms to represent the random profile with an accuracy of $10^{-2}$.
For example, when the period is $\Lambda=2\pi$, the root mean square is $\sigma=1/5$, and the correlation length is $l=1$, \cref{fig:lambda} shows the decay of the eigenvalues.

Hence, we can truncate the KL expansion and further reconstruct the profile of the random surface only considering the first several terms. Since the eigenfunctions $\{\varphi_j\}_{j\in \mathbb{N}}$ of the covariance operator are orthonormal and $\xi_0, \xi_{j,s}, \xi_{j,c}$ are independent and identically distributed random variables with zero mean and unit covariance, the eigenvalues of the covariance function are just that of the corresponding covariance matrix $C=(c_{ij})$ where $$c_{ij}=\mathbb{E}\left(\left\langle f-\mathbb{E} (f), \varphi_{i}\right\rangle\left\langle f-\mathbb{E} (f), \varphi_{j}\right\rangle\right).$$
Therefore, the second step of our reconstruction method is given as follows (\cref{alg:MCCUQ}). Note that the grating profile of a random periodic structure is reconstructed in the form of the sum of a finite Fourier series expansion and a truncated KL expansion represented by the reconstructed eigenvalues by this MCCUQ method.

\begin{algorithm}
	\caption{MCCUQ method}
	\label{alg:MCCUQ}
	\begin{algorithmic}
		\STATE{Define $\mathbf{c}=[\bar{c}_0, \bar{c}_1, \bar{c}_2,..., \bar{c}_{2k_{max}}]^{\top}$}
		\STATE{Set $\mathbf{c}=0$ (initializing)}
		\FOR{$m=1,2,...,M$ (sampling)}
		\STATE{Generate $\omega_m$, determine $y_0$ }
		\STATE{Define $\mathbf{c}_{m,k_{max}}=[c_{m,0},c_{m,1}...,c_{m,2k_{max}}]^{\mathrm{T}}$}
		\STATE{Set $c_{m,0}=y_0, c_{m,p}=0, p=1,2,...,2k_{max}$ (initializing)}
		\FOR{$k=1,2,...,k_{max}$}
		\STATE{Define $\mathbf{c}_{m,k}^{(0)}=[c_{m,0},c_{m,1}...,c_{m,2k}]^T$}
		\FOR{$t=0,1,2,...T$ (iteration)}
		\STATE{$\alpha=k sin(\theta_l)$,~$\beta=k cos(\theta_l)$,~$\theta_l\in (-\pi/2, \pi/2), l=1,2,...,L$}
		\STATE{Define
			$J_l(\mathbf{c}_{m,k}) =\left\|\displaystyle\sum_{n=-N}^{N} \psi_{n,l}(\omega_m) e^{i\alpha_n x + i\beta_n f(\omega_m;x)} +e^{i\alpha x-i\beta f(\omega_m,x)} \right\|_{\mathrm{L}^2(0,\Lambda)}^2$}
		\STATE{Define $\mathbf{J}(\mathbf{c}_{m,k}) = [ J_1(\mathbf{c}_{m,k}), ..., J_L(\mathbf{c}_{m,k}) ]^{\top}$}
		\STATE{Calculate the Jacobi matrix $\mathbf{DJ}=\left( \dfrac{\partial \mathbf{J}_l}{\partial c_{m,p}} \right)_{l=1,2,....,L,p=0,1,...,2k}$}
		\STATE{$\mathbf{c}_{m,k}^{(t+1)}=\mathbf{c}_{m,k}^{(t)}- \eta_k \mathbf{DJ}^{\top}(\mathbf{c}_{m,k}^{(t)}) \mathbf{J}(\mathbf{c}_{m,k}^{(t)})$}
		\ENDFOR
		\STATE{$c_{m,p}:=c_{m,p}^{(T+1)}, p=0,1,...,2k$}
		\ENDFOR
		\STATE{Let $ f_m(x)=c_{m,0} + \displaystyle\sum_{p=1}^{k_{max}} [c_{m,2p-1}\cos(px)+c_{m,2p}\sin(px)]$}
		\ENDFOR
		\STATE{$\mathbf{c} := \dfrac{1}{M}\displaystyle\sum_{m=1}^M \mathbf{c}_{m,k_{max}}$}
		\STATE{Let $\bar{f}(x)=\bar{c}_0 + \displaystyle\sum_{p=1}^{k_{max}} [\bar{c}_{2p-1}\cos(px)+\bar{c}_{2p}\sin(px)]$}
		\STATE{Calculate the covariance matrix $C$ where}
		\STATE{$c_{ij}:= \dfrac{1}{M}\left(\displaystyle\sum_{m=1}^M \left\langle f_m(x)-\bar{f}(x), \varphi_i(x)\right\rangle \left\langle f_m(x)-\bar{f}(x), \varphi_j(x)\right\rangle\right)$}
		\RETURN $\mathbf{\lambda}$ (eigenvalues of $C$)
	\end{algorithmic}
\end{algorithm}

\vspace{0.2cm}

\begin{remark}\label{RecoveryFormula}
Since the correlation length $l$ is greatly smaller than $\Lambda$, it can be deduced after simple calculation that the eigenvalues are only related to the root mean square $\sigma$ and the correlation length $l$ satisfying
$$\lambda_j \approx \sqrt{\pi} \sigma^2 l e^{-j^2 l^2/4},\ j=0,1,\ldots.$$
Thus we can further approximate these two statistics which characterize the random surface only by the eigenvalues reconstructed above. Specifically, we just need the first two eigenvalues $\lambda_0$ and $\lambda_1$, which leads to the following recovery formula
\begin{equation}
l'=\sqrt{4\ln \left(\dfrac{\lambda_0}{\lambda_1}\right)}, \quad \sigma' = \sqrt{\dfrac{\lambda_0}{\sqrt{\pi}l'}}.
\end{equation}
In fact, the correlation length and root mean square may be determined by any two eigenvalues.
\end{remark}

\section{Numerical results}
\label{sec:experiments}
To test the stability of our method, some random noise is added to all the measured scattering data, i.e., for each given sample, the scattering data takes the form
\begin{equation}
u(\omega;x,y_0):=u(\omega;x,y_0)(1+\tau \mbox{rand}),
\end{equation}
where rand represents uniformly distributed random numbers in $[-1,1]$ and $\tau$ is the noise level of the measured data. The near-field measurements $u(\omega;x,y_0)$ in this work are simulated by using an adaptive finite element (AFE) method with a perfectly matched absorbing layers\cite{Chen2003AdaptiveFEMPML,Bao2005Adaptive}. Since each sample is independent and identically distributed,  our numerical method can be performed in parallel.

For all of the numerical examples presented below, the number of eigenvalues determined by the covariance operator $C_f$ is chosen to satisfy the random surface with an accuracy of $10^{-4}$; the noise level $\tau$ and the regularization parameter $\gamma$ are taken as $0.1\%$ and $10^{-6}$ respectively;
the relaxation parameter dependent on the wavenumber $\kappa$ at each stage is chosen as $\eta_k=0.001/k^2$;
and the truncation of the infinite summation in Equation \cref{NonlinearProblem} is set to $N=8$. 

As for the number of realizations $M$,  as known to us, the Monte Carlo technique generally requires at least a sample size of $10^5$.
Besides, since the parameters we need to reconstruct are utilized to describe a random surface (a stationary Gaussian process), small sample size is insufficient for the stable reconstruction. Within the tolerance, we can appropriately reduce the sample size.  In our paper, since 1000 samples can reach our expectation, we set the number of realizations $M=1000$ in sampling the probability space.

\subsection{Accuracy of the algorithm}
\begin{figure}[htbp]
	\centering
	\includegraphics[width=10cm]{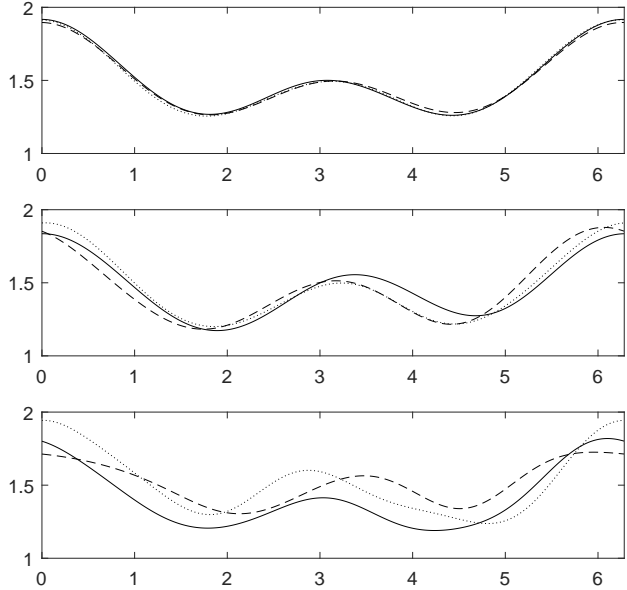}
	\caption{Certain realizations of randomly rough surface with $\sigma=1/15, 2/15, 1/5$(from the top to the bottom), respectively, for $l=1$.}\label{fig:surfaces}
\end{figure}

Recalling that the random surfaces are modeled by the KL expansion with the covariance function $$c(x-y)=\sigma^2 exp(-|x-y|^2/l^2),$$
to test the accuracy of our algorithm, we first show three realizations of randomly rough surfaces in one period  when $\tilde{f}(x)=1.5+0.2\cos(x)+0.2\cos(2x)$, the correlation length $l=1$, and the root mean square $\sigma=1/15, 2/15, 1/5$,
 respectively, as seen in
\cref{fig:surfaces}. It is clear that the deviation of the sample profile from the center profile becomes larger as the root mean square $\sigma$ increases.

\begin{figure}[htbp]
	\centering
	\includegraphics[width=13cm]{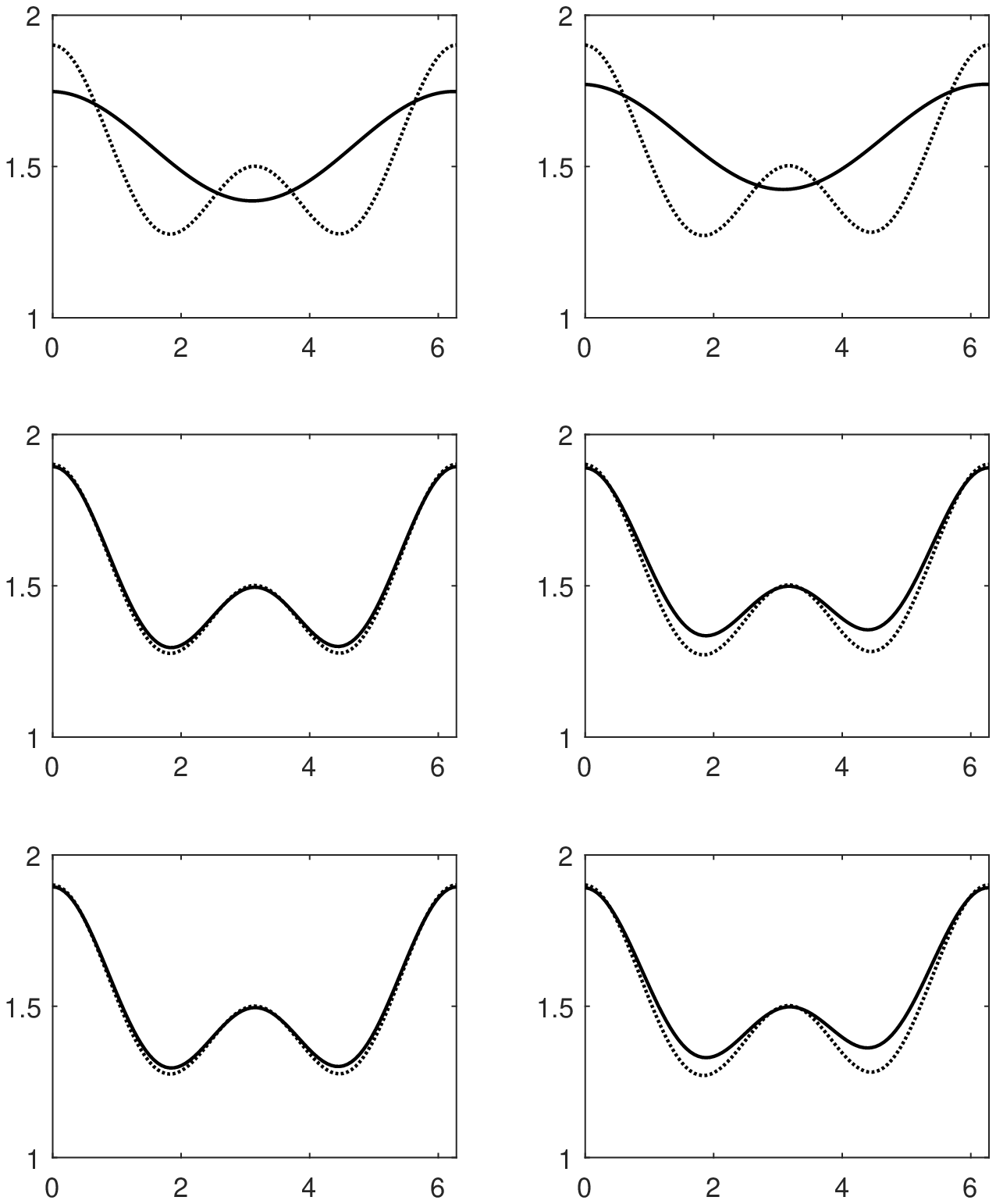}
	\caption{Evolutions of the reconstruction in \cref{example1}.  The solid curve: the reconstructed profile; the dotted curve: the given deterministic profile $\tilde{f}(x)$. From the top to the bottom are the reconstructions at $k=1$, $k=2$, and the reconstruction by our MCCUQ method. The left column is  $\sigma=1/15$; while the right column is $\sigma=1/5$.}\label{fig:ex1reconstruction}
\end{figure}

\begin{example}\label{example1}
To reconstruct a random periodic surface that is a sum of a deterministic function $$\tilde{f}(x)=1.5+0.2\cos(x)+0.2\cos(2x)$$ and a stationary Gaussian process with the covariance function $$c(x-y)={\sigma}^2 \exp(-|x-y|^2)$$
 with the correlation length $l=1$. Graphs of the original profile and the reconstructed profiles with different wavenumbers are shown in \cref{fig:ex1reconstruction}
with $\sigma=1/15$ and $1/5$, respectively. Since the original deterministic profile consists only of a few Fourier modes, only a small number of iterations are needed to get a good reconstruction both for $\sigma=1/15$ and $1/5$.
\end{example}

\begin{figure}[htbp]
	\centering
	\includegraphics[width=13cm]{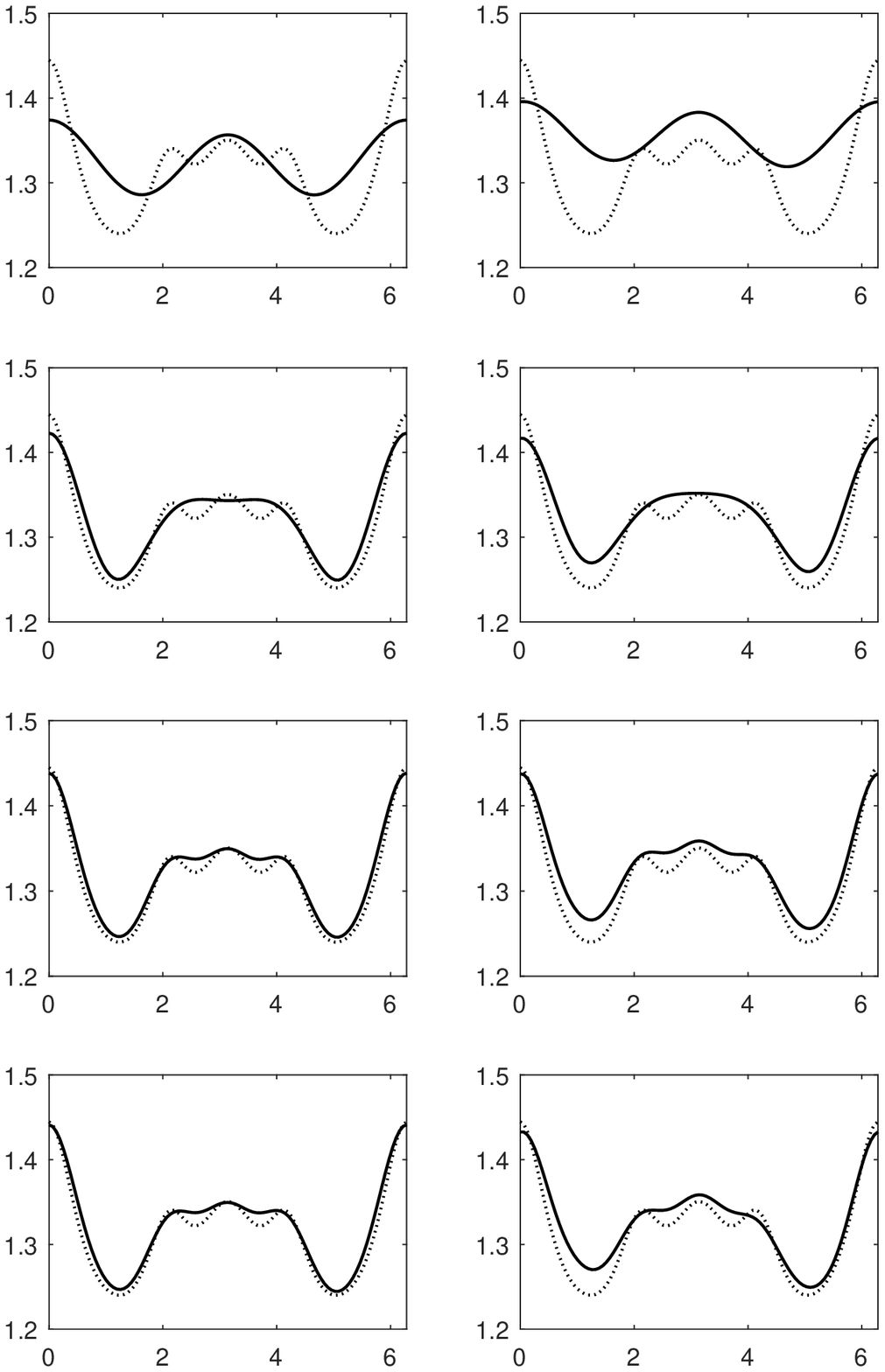}
	\caption{Evolutions of the reconstruction in  \cref{example2}.  The solid curve: the reconstructed profile; the dotted curve: the given deterministic profile $\tilde{f}(x)$. From the top to the bottom are reconstructions at $k=2$, $k=4$, $k=6$, and the reconstruction by our MCCUQ method, respectively. The left column is $\sigma=1/15$; while the right column is $\sigma=1/5$.}\label{fig:ex2reconstruction}
\end{figure}

\begin{example}\label{example2}
To reconstruct a random periodic surface that is a sum of a deterministic function $$\tilde{f}(x)=1.2+0.05\exp(\cos(2x))+0.04\exp(\cos(3x))$$ and a stationary Gaussian process with the covariance function $$c(x-y)={\sigma}^2 \exp(-|x-y|^2),$$
again with the correlation length  $l=1$.
Graphs of  the original profile and the reconstructed profiles with different wavenumbers are shown in \cref{fig:ex2reconstruction}.
As the original deterministic profile contains more Fourier modes, incident waves with higher frequencies are needed to get a better resolution of the reconstruction.
\end{example}

Recalling the formula of the expectation and standard deviation given in \cref{subsec:MC}, we denote $e(\omega)$ as the deviation of the reconstruction profile and the original deterministic profile for each sample. The mean error and standard deviation are denoted by $\bar{e}$ and $\bar{s}$ respectively.
In \cref{example1}, the mean error and standard deviation are $\bar{e}=1.0315\times 10^{-2}, \bar{s}=0.0639$ when the root mean square $\sigma=0.0667$, and $\bar{e}=3.5445\times 10^{-2}, \bar{s}=0.1778$ when the root mean square $\sigma=0.2$.
In \cref{example2}, $\bar{e}=7.9498\times 10^{-3}, \bar{s}=0.0654$ when the root mean square $\sigma=0.0667$, and $\bar{e}=1.4651\times 10^{-2}, \bar{s}=0.1951$ when the root mean square $\sigma=0.2$. Thus as shown in
\cref{fig:ex1reconstruction}-\cref{fig:ex2reconstruction}, all the reconstructions at least have an accuracy of $10^{-2}$  by our method (\cref{alg:MCC} and \cref{alg:MCCUQ}) proposed in this paper. Also, we can further evaluate the root mean square of the random structure directly by the standard deviation formula.

Therefore, we can reconstruct the unknown random periodic structure with an accuracy of $10^{-2}$ and give an evaluation of the root mean square value to certain degree, which shows a better reconstruction than that of the deterministic grating profile in \cite{bao2012ComputationalInverse}.

\subsection{Performance of the correlation length}
\begin{figure}[htbp]
	\centering
	\includegraphics[width=10cm]{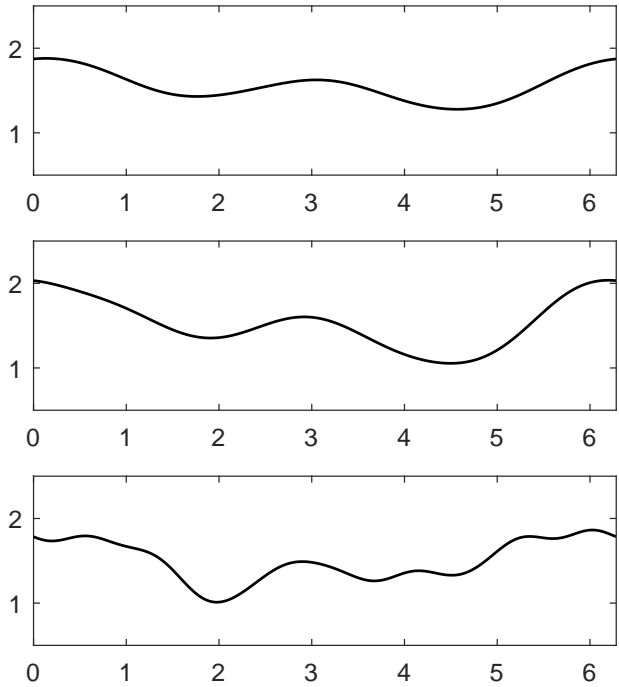}
	\caption{Certain realizations of the random surface with correlation lengths $l=1.5, 1.0, 0.5$(from the top to the bottom), respectively, for $\sigma=0.2$.}\label{CorrelationLength}
\end{figure}

We now discuss the influence of changes of the correlation length on our proposed method. For this, we consider random periodic surfaces with different correlation lengths. \cref{CorrelationLength} shows the realizations of random  periodic surfaces when $\tilde{f}(x)=1.5+0.2\cos(x)+0.2\cos(2x)$ and $l=1.5, 1.0, 0.5$ respectively. It is noticeable that the surface becomes rougher as the correlation length decreases.

We still consider the random periodic surfaces in \cref{example1} and \cref{example2} mentioned above.
By simulation, it is deduced that all the reconstructions of the random periodic surfaces for the correlation length $l=1.5, 1.0, 0.5$ and the root mean square $\sigma=1/15, 2/15, 1/5$ in \cref{example1} at the wavenumber $k=2$ and  in \cref{example2} at the wavenumber $k=6$ have an accuracy of $10^{-2}$.
That is, our proposed algorithm has the convergence and stability for both the rougher profile with $l=0.5$ and the larger root mean square with $\sigma=0.2$, which shows that the method proposed in this work is reliable and efficient.

\subsection{Evaluation of the statistics of the random periodic surface}

\begin{figure}[htbp]
	\centering
	\includegraphics[width=12cm]{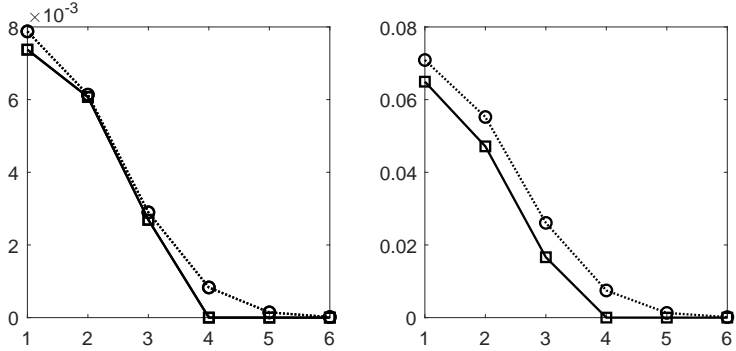}
	\caption{Reconstructions  of the eigenvalues in \cref{example1} by our MCCUQ method. The solid curve: the reconstructed eigenvalues; the dotted curve: the original eigenvalues. The left is $\sigma=1/15, k=2$; while the right is $\sigma=1/5, k=2$.}\label{ex1relambda}
\end{figure}

\begin{figure}[htbp]
	\centering
	\includegraphics[width=12cm]{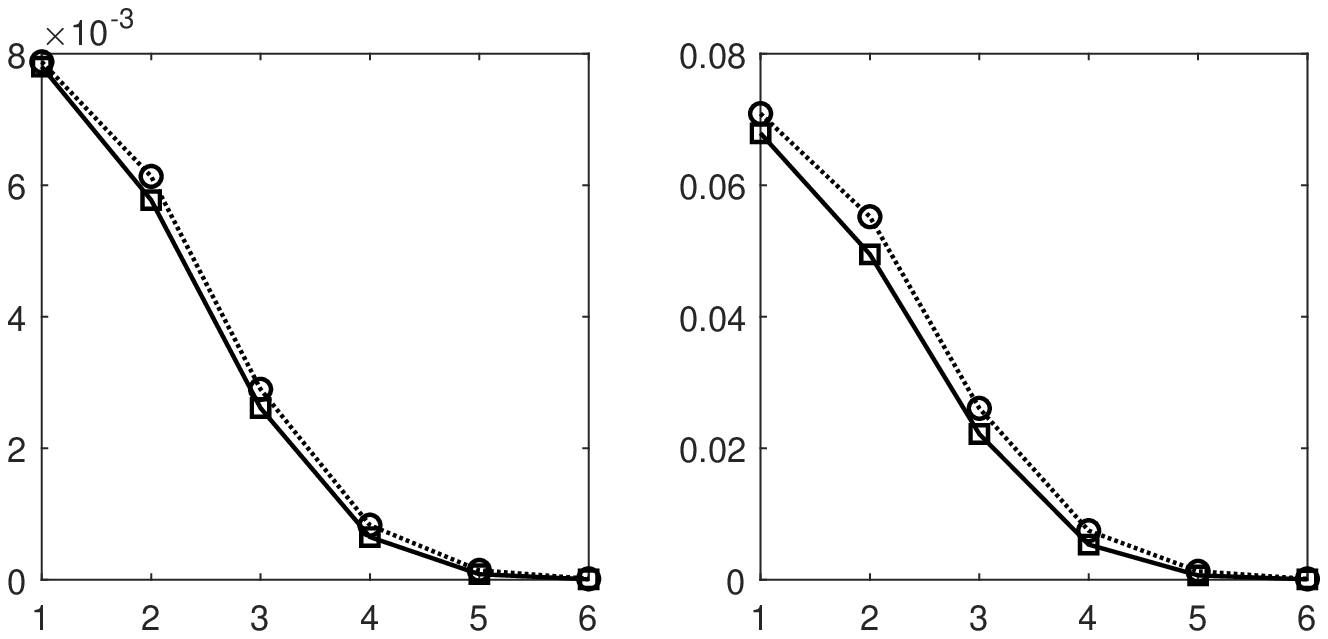}
	\caption{Reconstructions  of the eigenvalues in \cref{example2} by our MCCUQ method. Solid curve: the reconstructed eigenvalues; the dotted curve: the original eigenvalues. The left is $\sigma=1/15, k=6$ and the right is $\sigma=1/5, k=6$.}\label{ex2relambda}
\end{figure}

We still consider the random periodic surfaces in \cref{example1} and \cref{example2} mentioned above.  \cref{ex1relambda} shows the reconstruction of eigenvalues in \cref{example1} when $\sigma=1/15$ and $\sigma=1/5$ respectively.
 \cref{ex2relambda} shows the reconstruction of eigenvalues in \cref{example2} when $\sigma=1/15$ and $\sigma=1/5$ respectively.
As shown in \cref{ex1relambda}-\cref{ex2relambda}, all the eigenvalues reconstructed by our MCCUQ method achieve at least an accuracy of $10^{-3}$. By doing more experiments, we observe that with the smaller root mean square or the larger correlation length, the reconstructed eigenvalues can reach a higher accuracy.

From the reconstructed eigenvalues, we can calculate the statistics of the random periodic surface by the recovery formula discussed in \cref{RecoveryFormula}.
It can be calculated directly that $l'=1.0415$, $\sigma'=0.0657$ when $\sigma=0.0667$, and $l'=1.1321$, $\sigma'=0.1799$ when $\sigma=0.2$ in \cref{example1},
$l'=1.0282$, $\sigma'=0.0642$ when $\sigma=0.0667$, and $l'=0.9920$, $\sigma'=0.1955$ when $\sigma=0.2$ in \cref{example2}.
Since we know the original correlation length is $l=1$ in \cref{example1} and \cref{example2}, it is clear that the evaluation of the two key statistics of the random periodic surfaces by our method is excellent.

Numerical results show that the recovery formula proposed in \cref{RecoveryFormula} can reconstruct the root mean square and the correlation length of the random structure explicitly by the eigenvalues reconstructed by our MCCUQ method (\cref{alg:MCCUQ}).
Also, the conclusion is deduced that the larger the correlation length $l$ is or the smaller the root mean square $\sigma$ is, the more accurate the evaluation of the statistics is. This means the smoother the random surfaces are, the better the evaluations of the statistics are.

\section{Concluding remarks}
\label{sec:conclusions}
We have presented an efficient numerical method for solving the inverse scattering problem by random periodic structure. Our method (MCCUQ) is based on a novel combination of the Monte Carlo technique for sampling the probability space, a continuation method with respect to the wavenumber, and the Karhunen-Lo$\grave{e}$ve expansion of the random structure, which reconstructs key statistical properties of the profile for the unknown random periodic structure from boundary measurements of the scattered fields away from the structure. Numerical results are presented to demonstrate the reliability and efficiency of the proposed method.
For the direct problem in this work, we use the finite element method to find the scattered field. Since the recovery formula is explicit, the majority of the computing time in our simulation is spent  on the calculation of the direct problem  and we can calculate the statistics quickly with the observation data.
 To the best of knowledge, this is the first attempt to directly solve the inverse scattering problem to determine a random interface or obstacle. There are many interesting open problems along this direction.   Although our numerical examples demonstrate the convergence of MCCUQ, no rigorous convergence and stability analysis is available at present. Obviously new theoretical framework must be developed. Another important future project is to develop algorithms for solving the inverse scattering problems by penetrable randomly rough surfaces. Also, more realistic inverse scattering for three-dimensional Maxwell's equations should also be investigated.
In the three dimensional problem, the finite element method for the direct problem will increase exponentially the computational cost while the Monte Carlo technique and the inversion procedure are not affected too much by the dimensionality. Hence it is promising to extend our reconstruction method to higher dimensional problems.
Finally, it would be also significant to develop numerical methods for solving inverse scattering by non-periodic randomly rough surfaces.

\bibliographystyle{siamplain}
\bibliography{references}

\end{document}